\documentclass[11pt]{article}
\usepackage[brazil,english]{babel}
\usepackage[latin1]{inputenc}
\usepackage{fullpage}
\usepackage[dvips]{graphicx}
\usepackage{amssymb}
\usepackage{amsmath}
\usepackage{color}
\usepackage{footnote}
\usepackage{hyperref}

\newtheorem{theorem}{Theorem}[section] 
  
\newtheorem{proposition}[theorem]{Proposition}  
\newtheorem{lemma}[theorem]{Lemma}  
\newtheorem{definition}[theorem]{Definition} 
\newtheorem{example}[theorem]{Example}  
\newtheorem{remark}[theorem]{Remark}  

\def\beginproof{\noindent{\em Proof}.~}
\def\endproof{{\ \hfill\hbox{\fbox{}}\parfillskip 0pt}\par}

\renewcommand{\thefootnote}{\arabic{ftnote}}

\newcommand{\R}{\mathbb{R}}
\definecolor{cmg}{RGB}{250, 0, 250}

\title{A sequential optimality condition for Mathematical Programs with Cardinality 
Constraints}

\author{Evelin H. M. Krulikovski\footnotemark[3] \and Ademir A. Ribeiro\footnotemark[4]
\footnotemark[5] \and Mael Sachine\footnotemark[4]}

\begin{document}

\maketitle
\renewcommand{\thefootnote}{\fnsymbol{footnote}}
\footnotetext[3]{Graduate Program in Mathematics, Federal University of Paran\'a, 
Brazil (evelin.hmk@gmail.com).}
\footnotetext[4]{Department of Mathematics, Federal University of Paran\'a, 
Brazil (ademir.ribeiro@ufpr.br, mael@ufpr.br).}
\footnotetext[5]{Partially supported by CNPq, Brazil, Grant 309437/2016-4.}

\renewcommand{\thefootnote}{\arabic{ftnote}}

\begin{abstract}
In this paper we propose an Approximate Weak stationarity ($AW$-stationarity) concept 
designed to deal with {\em Mathematical Programs with Cardinality Constraints} (MPCaC), 
and we proved that it is a legitimate optimality condition independently of any constraint 
qualification. Such a sequential optimality condition improves weaker stationarity conditions, 
presented in a previous work. Many research on sequential optimality conditions has been 
addressed for nonlinear constrained optimization in the last few years, some works in the 
context of MPCC and, as far as we know, no sequential optimality condition has been proposed 
for MPCaC problems. We also establish some relationships between our $AW$-stationarity 
and other usual sequential optimality conditions, such as AKKT, CAKKT and PAKKT.
We point out that, despite the computational appeal of the sequential optimality 
conditions, in this work we are not concerned with algorithmic consequences. 
Our aim is purely to discuss theoretical aspects of such conditions for MPCaC problems.
\end{abstract}


{\bf Subclass. }
90C30, 90C33, 90C46 

\thispagestyle{plain}

\section{Introduction}
In this paper we propose a sequential optimality condition, associated to the weak 
stationarity condition presented in our previous work \cite{KrulikovskiRibeiroSachine20aX}, 
designed to deal with {\em Mathematical Programs with Cardinality Constraints} (MPCaC) 
given by 
\begin{equation}
\label{aw_prob:mpcac}
\begin{array}{cl}
\displaystyle\mathop{\rm minimize }  & f(x)  \\
{\rm subject\ to } & x \in X, \\
& \|x\|_0\leq \alpha,
\end{array}
\end{equation}
where $f:\R^n\to\R$ is a continuously differentiable function, $X\subset\R^n$ is a set given 
by equality and/or inequality constraints, $\alpha>0$ is a given natural number and 
$\|x\|_0$ denotes the cardinality of the vector $x\in\R^n$, that is, the number of nonzero 
components of $x$. We assume that $\alpha<n$, since otherwise the cardinality constraint 
would be innocuous. Note, however, that if $\alpha$ is too small, the cardinality constraint 
may be too restrictive leading to an empty feasible set.
Furthermore, the main difference between the problem (\ref{aw_prob:mpcac}) and a standard nonlinear 
programming problem is that the cardinality constraint, despite of the notation, is not a norm, 
nor continuous neither convex. 

One reformulation to deal with this difficult cardinality constraint consists of addressing 
its continuous counterpart \cite{BurdakovKanzowSchwartz16}
\begin{equation}
\label{aw_prob:relax}
\begin{array}{cl}
\displaystyle\mathop{\rm minimize }_{x,y} & f(x)  \\
{\rm subject\ to } & x \in X, \\
& e^Ty\geq n-\alpha, \\
& x_iy_i=0,\; i=1,\ldots,n, \\
& 0\leq y_i \leq 1, \; i=1,\ldots,n,
\end{array}
\end{equation}
which will be referred to as relaxed problem and, with some abuse of terminology, 
will be indicated as MPCaC as well. 
It can be seen that these problems are equivalent in the sense that global 
solutions of (\ref{aw_prob:mpcac}) correspond, in a natural way, to global solutions 
of (\ref{aw_prob:relax}) and, if $x^*\in\R^n$ is a local minimizer of (\ref{aw_prob:mpcac}), 
then every feasible pair $(x^*,y^*)$ is a local minimizer of (\ref{aw_prob:relax}).

In \cite{KrulikovskiRibeiroSachine20aX} we proposed new and weaker stationarity 
conditions for this class of problems, by means of a unified approach that goes 
from the weakest to the strongest stationarity. Indeed we cannot assert about 
KKT points for MPCaC problems, since some standard constraint qualifications 
are violated. This occurs in view of the complementarity constraints $x_iy_i=0$, 
$i=1,\ldots,n$. 
However, the weaker condition proposed in \cite{KrulikovskiRibeiroSachine20aX}, 
called $W_{I}$-stationarity, even being weaker than KKT, is not a necessary 
optimality condition. Therefore, we 
propose in this work an Approximate Weak stationarity ($AW$-stationarity) concept, 
which will be proved to be a legitimate optimality condition, independently of 
any constraint qualification. 

In the last few years, special attention has been paid to the so-called sequential 
optimality conditions for nonlinear constrained optimization 
\cite{AndreaniFazzioSchuverdtSecchin,AndreaniHaeserMartinez,AndreaniMartinezRamosSilva16,AndreaniMartinezRamosSilva18,AndreaniMartinezSvaiter,MartinezSvaiter,RibeiroSachineSantos18}. 
Sequential optimality conditions are intrinsically related to the stopping criteria 
of numerical algorithms, and their study aims at unifying the theoretical convergence 
analysis associated with the corresponding algorithm. Within this context, for 
instance, the augmented Lagrangian method (see \cite{BirginMartinez} and references therein) 
has been extensively analyzed, being shown to satisfy weak sequential conditions, 
thus giving rise to strong convergence results.

Sequential optimality conditions are necessary for optimality, i.e., a local 
minimizer of the problem under consideration verifies such a condition, independently 
of the fullfilment of any constraint qualification (CQ). 
The \emph{approximate Karush-Kuhn-Tucker} (AKKT) is one of the most popular of these 
conditions, and it was defined in \cite{AndreaniHaeserMartinez} 
and \cite{QiWei}. Another two sequential optimality conditions for standard nonlinear 
programming, both stronger than AKKT, are 
\emph{positive approximate KKT} (PAKKT) \cite{AndreaniFazzioSchuverdtSecchin}
and \emph{complementary approximate KKT} (CAKKT) \cite{AndreaniMartinezSvaiter}. 
Whenever it is proved that an AKKT (or CAKKT or PAKKT) point is indeed a 
Karush-Kuhn-Tucker (KKT) point under a certain CQ, 
any algorithm that reaches AKKT (or CAKKT or PAKKT) points 
(e.g. augmented Lagrangian-type methods) automatically has the theoretical convergence 
established assuming the same CQ. This paves the grounds for the aforementioned unification. 

Sequential optimality conditions have also been proposed for nonstandard optimization 
\cite{AndreaniHaeserSecchinSilva,HelouSantosSimoes20a,HelouSantosSimoes20b,Ramos}. 
In the context of {\em Mathematical Programs with Equilibrium Constraints} (MPECs) 
and motivated by AKKT, it was introduced in \cite{Ramos} the MPEC-AKKT condition with a 
geometric appeal and in \cite{AndreaniHaeserSecchinSilva}, new conditions were established 
for {\em Mathematical Problems with Complementarity Constraints} (MPCCs), 
namely $AW$-, $AC$- and $AM$-stationarity. The latter one was compared with the sequential 
condition present in \cite{Ramos}.

Even though there is a considerable literature devoted to sequential conditions for 
standard nonlinear optimization and even for specific problems (MPCC and MPEC), 
to the best of our knowledge, no sequential optimality condition has been 
proposed for MPCaC problems. Such problems are very degenerate because of the problematic 
complementarity constraints $x_iy_i=0$ and therefore the known sequential optimality 
conditions may not be suitable to deal with them. Thereby, we propose a sequential 
optimality condition, namely $AW$-stationarity, associated to $W_{I}$-stationarity 
and designed to deal with MPCaC problems. This condition is based on the one 
proposed in \cite{AndreaniHaeserSecchinSilva} for MPCC problems.  
The main contribution of this paper is that $AW$-stationarity is indeed a necessary 
optimality condition, without any constraint qualification assumption. 
We also establish some relationships between our $AW$-stationarity and other well known 
sequential optimality conditions. In particular, and surprisingly, we prove that 
AKKT fails to detect good candidates for optimality  
for every MPCaC problem.

We stress that, despite the algorithmic appeal of the sequential optimality 
conditions, in this work we are neither concerned with
applications nor with computational aspects or algorithmic consequences. 
Our aim is to discuss theoretical aspects of such conditions for MPCaC problems.

The paper is organized as follows: in Section \ref{aw_sec:prelim} we establish the 
notation, some definitions and results concerning standard nonlinear programming 
and recall the weak stationarity concept proposed in our previous 
work \cite{KrulikovskiRibeiroSachine20aX}. Section \ref{aw_sec:aw} presents 
the main results of this paper, concerning sequential optimality conditions for MPCaC. 
In Section \ref{aw_sec:rel_akkt_cakkt} we provide some relationships between approximate 
stationarity for standard nonlinear optimization and $AW$-stationarity. 
Concluding remarks are presented in Section~\ref{aw_sec:concl}.

\medskip

{\noindent\bf Notation.} Throughout this paper, for vectors $x,y\in\R^n$, $x*y$ denotes 
the Hadamard product between $x$ and $y$, that is, the vector obtained by the 
componentwise product of $x$ and $y$. In the same way, the ``min'' in the vector 
$\min\{x,y\}$ is taken componentwise. We also use the following sets of indices: 
$I_{00}(x,y)=\{i\mid x_i=0,y_i=0\}$, 
$I_{\pm 0}(x,y)=\{i\mid x_i\neq 0,y_i=0\}$,
$I_{0+}(x,y)=\{i\mid x_i=0,y_i\in(0,1)\}$,
$I_{01}(x,y)=\{i\mid x_i=0,y_i=1\}$,
$I_{0\,>}(x,y)=\{i\mid x_i=0,y_i>0\}$ and 
$I_0(x)=\{i\mid x_i=0\}$. For a vector-valued function $\xi:\R^n\to\R^s$, denote 
$I_\xi(x)=\{i\mid \xi_i(x)=0\}$, the set of active indices, and 
$\nabla\xi= (\nabla\xi_1\ldots\nabla\xi_s)$, the transpose of the Jacobian of $\xi$.

\section{Preliminaries}
\label{aw_sec:prelim}
In this section we recall some basic definitions and results related to standard 
nonlinear programming (NLP), as well as the weak stationarity concept proposed in 
our previous work \cite{KrulikovskiRibeiroSachine20aX}. 

Consider first the problem 
\begin{equation}
\label{aw_prob:nlp}
\begin{array}{cl}
\displaystyle\mathop{\rm minimize } & f(x)\\
{\rm subject\ to } &  g(x)\leq 0, \\
& h(x)=0,
\end{array}
\end{equation}
where $f:\R^n\to\R$, $g:\R^n\to\R^m$ and $h:\R^n\to\R^p$ are continuously differentiable 
functions. The feasible set of the problem (\ref{aw_prob:nlp}) is denoted by 
\begin{equation}
\label{aw_feas_set}
\Omega=\{x\in\R^n\mid g(x)\leq 0, h(x)=0\}.
\end{equation}
\begin{definition}
We say that $x^*\in\Omega$ is a global solution of the problem (\ref{aw_prob:nlp}), that is, 
a global minimizer of $f$ in $\Omega$, when $f(x^*)\leq f(x)$ for all $x\in\Omega$. If 
$f(x^*)\leq f(x)$ for all $x\in\Omega$ such that $\|x-x^*\|\leq\delta$, for some constant 
$\delta>0$, $x^*$ is said to be a local solution of the problem.
\end{definition}

A feasible point $x^*\in\Omega$ is said to be {\em stationary} for the problem (\ref{aw_prob:nlp}) 
if there exists a vector $\lambda=(\lambda^g,\lambda^h)\in\R_+^m\times\R^p$ 
(Lagrange multipliers) such that 
\begin{subequations}
\begin{align}
\nabla f(x^*)+\sum_{i=1}^{m}\lambda_i^g\nabla g_i(x^*)
+\sum_{i=1}^{p}\lambda_i^h\nabla h_i(x^*)=0, & \label{aw_kkt_grad} \\ 
(\lambda^g)^Tg(x^*)=0. \label{aw_kkt_compl}
\end{align}
\end{subequations}
The function $L:\R^n\times\R^m\times\R^p\to\R$ given by 
\begin{equation}
\label{aw_lagr:nlp}
L(x,\lambda^g,\lambda^h)=f(x)+(\lambda^g)^Tg(x)+(\lambda^h)^Th(x)
\end{equation}
is the \textit{Lagrangian} function associated with the problem (\ref{aw_prob:nlp}).

The conditions (\ref{aw_kkt_grad})--(\ref{aw_kkt_compl}) are known as Karush-Kuhn-Tucker (KKT) 
conditions and, under certain qualification assumptions, are satisfied at a minimizer.

\subsection{Constraint qualifications}
\label{aw_sec:cq}
There are a lot of constraint qualifications, that is, conditions under which every minimizer 
satifies KKT. In order to discuss some of them, let us recall the definition of cone, which 
plays an important role in this context. 

We say that a nonempty set $C\subset\R^n$ is a {\em cone} if $td\in{C}$ for all $t\geq 0$ 
and $d\in{C}$. Given a set $S\subset\R^n$, its {\em polar} is the cone 
$$
S^\circ=\{p\in\R^n\mid p^Tx\leq 0,\ \forall x\in S\}.
$$

Associated with the feasible set of the problem (\ref{aw_prob:nlp}), we have 
the {\em tangent cone} 
$$
T_{\Omega}(\bar{x})=\left\{d\in\R^n\mid\exists(x^k)\subset\Omega\mbox{, } 
(t_k)\subset\R_+ \mbox{ : }t_k\to 0 \mbox{ and } 
\dfrac{x^k-\bar{x}}{t_k}\to d\right\}
$$
and the {\em linearized cone} 
$$
D_{\Omega}(\bar{x})=\left\{d\in\R^n\mid\nabla g_i(\bar{x})^Td \leq  0, \;
i\in I_g(\bar{x})
\mbox{ and }\nabla h(\bar{x})^Td=0 \right\}.
$$

The following basic result says that we may ignore inactive constraints when 
dealing with the tangent and linearized cones. 
\begin{lemma}
\label{aw_lm:inactive}
Consider a feasible point $\bar{x}\in\Omega$, an index set $J\supset I_g(\bar{x})$ 
and 
$$
\Omega'=\{x\in\R^n\mid g_{i}(x)\leq 0,\, i\in J, \; h(x)=0\}.
$$
Then, $T_{\Omega}(\bar{x})=T_{\Omega'}(\bar{x})$ and 
$D_{\Omega}(\bar{x})=D_{\Omega'}(\bar{x})$.
\end{lemma}
\beginproof
Note first that $\bar{x}\in\Omega'$ since $\Omega\subset\Omega'$. Moreover, since 
$g_i(\bar{x})<0$ for $i\notin J$, there exists $\delta>0$ such that 
$B(\bar{x},\delta)\cap\Omega'=B(\bar{x},\delta)\cap\Omega$. Thus, 
$T_{\Omega'}(\bar{x})=T_{\Omega}(\bar{x})$ because the conditions $t_k\to 0$ and 
$(x^k-\bar{x})/t_k\to d$ imply that $x^k\to\bar{x}$. 
The equality between the linearized cones is straightforward, as the active 
indices corresponding to $\Omega$ and $\Omega'$ coincide.
\endproof

Now we relate the cones of feasible sets when some variables do not appear in the constraints. 

\begin{lemma}
\label{aw_lm:withoutx}
Consider the general feasible set $\Omega$, defined in (\ref{aw_feas_set}), and the set 
$$
\Omega'=\{(x,y)\in\R^n\times\R^m\mid g(x)\leq 0, h(x)=0\}.
$$
Given a feasible point $(\bar{x},\bar{y})\in\Omega'$, we have 
$$
T_{\Omega'}(\bar{x},\bar{y})=T_{\Omega}(\bar{x})\times\R^m
\quad\mbox{and}\quad 
D_{\Omega'}(\bar{x},\bar{y})=D_{\Omega}(\bar{x})\times\R^m.
$$
As a consequence, 
$$
T_{\Omega'}^\circ(\bar{x},\bar{y})=T_{\Omega}^\circ(\bar{x})\times\{0\}
\quad\mbox{and}\quad 
D_{\Omega'}^\circ(\bar{x},\bar{y})=D_{\Omega}^\circ(\bar{x})\times\{0\}.
$$
\end{lemma}
\beginproof
The relation between the tangent cones follows directly from the definition. 
Moreover, if $\zeta(x,y)=g(x)$ and $\xi(x,y)=h(x)$ represent the constraints that 
define $\Omega'$ and $d=(\alpha,\beta)$, then 
$$
\nabla\zeta_i(x,y)^Td=\nabla g_i(x)^T\alpha\quad\mbox{and}\quad 
\nabla\xi_j(x,y)^Td=\nabla h_j(x)^T\alpha,
$$
which gives the second claim. Finally, the last statement of the lemma follows from 
the fact that $(S\times\R^m)^\circ=S^\circ\times\{0\}$ for any set $S\subset\R^n$.
\endproof

The two weakest constraint qualifications are defined below. 
\begin{definition}
\label{aw_acq_gcq}
We say that Abadie constraint qualification (ACQ) holds at $\bar{x}\in\Omega$ if 
$T_{\Omega}(\bar{x})=D_{\Omega}(\bar{x})$. If 
$T_{\Omega}^\circ(\bar{x})=D_{\Omega}^\circ(\bar{x})$, we say that Guignard constraint 
qualification (GCQ) holds at $\bar{x}$.
\end{definition}

In the following lemma we analyze GCQ for simple complementarity constraints. 

\begin{lemma}
\label{aw_lm:compl1}
Consider the set 
$$
\Omega=\{(x,y)\in\R^n\times\R^n\mid y\geq 0, x*y=0\}.
$$
Given $(\bar{x},\bar{y})\in\Omega$, there holds 
$T_{\Omega}^\circ(\bar{x},\bar{y})=D_{\Omega}^\circ(\bar{x},\bar{y})$.
\end{lemma}
\beginproof
Denote the constraints that define $\Omega$ by 
$\zeta(x,y)=-y$ and $\xi(x,y)=x*y$. Given $d=(u,v)\in D_{\Omega}(\bar{x},\bar{y})$, 
we claim that the vectors $d_1=(u,0)$ and $d_2=(0,v)$ belong to 
$T_{\Omega}(\bar{x},\bar{y})$. Indeed, since 
$$
\bar{y}_iu_i+\bar{x}_iv_i=
\nabla\xi_i(\bar{x},\bar{y})^Td=0
$$
for all $i=1,\ldots,n$, we have $u_{I_{0>}}=0$ and $v_{I_{\pm0}}=0$, where we used 
the simplified notation $I_{0>}=I_{0>}(\bar{x},\bar{y})$ and  
$I_{\pm 0}=I_{\pm 0}(\bar{x},\bar{y})$. Thus, the sequences 
$t_k={1}/{k}$ and $(x^k,y^k)=(\bar{x}+t_ku,\bar{y})$ satisfy $y^k\geq 0$, $x^k_{I_{0>}}=0$, 
$y^k_{I_{\pm0}\cup I_{00}}=0$, which means that $(x^k,y^k)\subset\Omega$, and 
$$
\dfrac{(x^k,y^k)-(\bar{x},\bar{y})}{t_k}\to{d_1},
$$ 
proving that $d_1\in T_{\Omega}(\bar{x},\bar{y})$. 
Now, defining $(z^k,w^k)=(\bar{x},\bar{y}+t_kv)$, we have 
$$
\dfrac{(z^k,w^k)-(\bar{x},\bar{y})}{t_k}\to{d_2}\quad\mbox{and}\quad z^k*w^k=0.
$$
Furthermore, for $i\in I_{0>}$ we have $w_i^k>0$ for all sufficiently large $k$. 
On the other hand, if $i\in I_{00}$, then the constraint $\zeta_i$ is active and hence, 
$
-v_i=\nabla\zeta_i(\bar{x},\bar{y})^Td\leq 0,
$
giving $w_i^k=\bar{y}_i+t_kv_i=t_kv_i\geq 0$. Thus, $(z^k,w^k)\subset\Omega$, 
which yields ${d_2}\in T_{\Omega}(\bar{x},\bar{y})$, proving the claim. 
Finally, given $p\in T_{\Omega}^\circ(\bar{x},\bar{y})$ we conclude 
that $p^Td=p^Td_1+p^Td_2\leq 0$, proving that $p\in D_{\Omega}^\circ(\bar{x},\bar{y})$.
\endproof

\subsection{Sequential optimality conditions for standard NLP}
\label{aw_sec:akkt}
The goal of this section is to present some well known approximate optimality conditions 
for nonlinear constrained optimization 
\cite{AndreaniHaeserMartinez,AndreaniMartinezRamosSilva16,AndreaniMartinezRamosSilva18,AndreaniMartinezSvaiter,BirginMartinez,MartinezSvaiter}. 

\begin{definition}
\label{aw_def:akkt}
Let $\bar{x}\in\R^n$ be a feasible point for the problem (\ref{aw_prob:nlp}). We say that 
$\bar{x}$ is an \textnormal{Approximate KKT} (AKKT) point if 
there exist sequences $(x^k)\subset{\mathbb R}^n$ and 
$(\lambda^k)=\big(\lambda^{g,k},\lambda^{h,k}\big)\subset\R_+^m\times\R^p$ 
such that $x^k\to \bar{x}$, 
\begin{subequations}
\begin{align}
\nabla_xL(x^k,\lambda^{g,k},\lambda^{h,k})\to 0, \label{aw_grad1xxx} \\ 
\min\{-g(x^k),\lambda^{g,k}\}\to 0.  \label{aw_compl1xxx} 
\end{align}
\end{subequations}
\end{definition}

We have below two stronger conditions than AKKT.
\begin{definition}
\label{aw_def:cakkt}
Let $\bar{x}\in\R^n$ be a feasible point for the problem (\ref{aw_prob:nlp}). We say that 
$\bar{x}$ is a \textnormal{Complementary Approximate KKT} (CAKKT) point if there 
exist sequences $(x^k)\subset{\mathbb R}^n$ and 
$(\lambda^k)=\big(\lambda^{g,k},\lambda^{h,k}\big)\subset\R_+^m\times\R^p$ such 
that $x^k\to \bar{x}$, 
\begin{subequations}
\begin{align}
\nabla_xL(x^k,\lambda^{g,k},\lambda^{h,k})\to 0, 
\label{aw_gradcakkt} \\ \lambda^{g,k}*g(x^k)\to 0\quad\mbox{and}\quad
\lambda^{h,k}*h(x^k)\to 0. \label{aw_complcakkt} 
\end{align}
\end{subequations}
\end{definition}

\begin{remark}
\label{aw_rm:compl4_compl1}
Note that if $(\alpha^k)\subset{\mathbb R_+}$ and $(\beta^k)\subset{\mathbb R}$ are 
sequences satisfying $\alpha^k\beta^k\to 0$ and $\beta^k\to\bar\beta\leq 0$, then 
$\min\{-\beta^k,\alpha^k\}\to 0$. Indeed, if $\bar\beta<0$, we have $\alpha^k\to 0$ and 
hence $\alpha^k<-\beta^k$ for all $k$ sufficiently large, giving
$\min\{-\beta^k,\alpha^k\}=\alpha^k\to 0$. On the other side, if $\bar\beta=0$, we also 
conclude that $\min\{-\beta^k,\alpha^k\}\to 0$, since $\alpha^k\geq 0$. This means that 
condition (\ref{aw_complcakkt}) implies (\ref{aw_compl1xxx}), and thus CAKKT implies AKKT.
\end{remark}

Another known sequential optimality condition relates to the sign of the multipliers. 
\begin{definition}
\label{aw_def:pakkt}
Let $\bar{x}\in\R^n$ be a feasible point for the problem (\ref{aw_prob:nlp}). We say that 
$\bar{x}$ is a \textnormal{Positive Approximate KKT} (PAKKT) point if there 
exist sequences $(x^k)\subset{\mathbb R}^n$ and 
$(\lambda^k)=\big(\lambda^{g,k},\lambda^{h,k}\big)\subset\R_+^m\times\R^p$ such 
that $x^k\to \bar{x}$, 
\begin{subequations}
\begin{align}
\nabla_xL(x^k,\lambda^{g,k},\lambda^{h,k})\to 0, 
\label{aw_gradpakkt} \\ 
\min\{-g(x^k),\lambda^{g,k}\}\to 0, \label{aw_complpakkt} \\ 
\lambda_i^{g,k}g_i(x^k)>0 \mbox{ if } \displaystyle\mathop{\rm lim\, sup}_{k\to\infty}
\frac{\lambda_i^{g,k}}{\delta_k}>0, \label{aw_pospakkt1}  \\ 
\lambda_j^{h,k}h_j(x^k)>0 \mbox{ if } \displaystyle\mathop{\rm lim\, sup}_{k\to\infty}
\frac{|\lambda_j^{h,k}|}{\delta_k}>0, \label{aw_pospakkt2} 
\end{align}
\end{subequations}
where $\delta_k= \|(1,\lambda^k)\|_\infty$. 
\end{definition}

As well known in the literature, all the three sequential conditions above are 
necessary optimality conditions without any constraint qualification.

\subsection{Weak stationarity for MPCaC}
\label{aw_sec:sequential}
In this section we recall the weaker stationarity concept and some related results, 
established in \cite{KrulikovskiRibeiroSachine20aX}, for MPCaC. As we have seen in that work, 
except in special cases, e.g., when $X$ is given by linear constraints, we do not have the 
fulfillment of constraint qualifications for the relaxed problem (\ref{aw_prob:relax}). 
So, the standard KKT conditions are not necessary optimality conditions, fact that in turn 
justifies the study of weaker conditions.

For ease of presentation consider the functions $\theta:\R^n\to\R$, 
$G,H,\tilde{H}:\R^n\to\R^n$ given by 
$$
\theta(y)=n-\alpha-e^Ty\,,\quad G(x)=x\,,\quad H(y)=-y\quad\mbox{and}\quad\tilde{H}(y)=y-e. 
$$
Then we can rewrite the relaxed problem (\ref{aw_prob:relax}) as 
\begin{equation}
\label{aw_prob:relax1}
\begin{array}{cl}
\displaystyle\mathop{\rm minimize }_{x,y} & f(x)  \\
{\rm subject\ to } & g(x)\leq 0, h(x)=0, \\
& \theta(y)\leq 0, \\
& H(y)\leq 0, \tilde{H}(y)\leq 0, \\
& G(x)*H(y)=0.
\end{array}
\end{equation}

Given a feasible point $(\bar{x},\bar{y})$ for the problem (\ref{aw_prob:relax1}) and a set of 
indices $I$ such that 
\begin{equation}
\label{aw_indexI}
I_{0+}(\bar{x},\bar{y})\cup I_{01}(\bar{x},\bar{y})\subset I\subset I_0(\bar{x}),
\end{equation}
we have that $i\in I$ or $i\in I_{00}(\bar{x},\bar{y})\cup I_{\pm 0}(\bar{x},\bar{y})$
for all $i\in\{1,\ldots,n\}$. Thus, $G_i(\bar{x})=0$ or $H_i(\bar{y})=0$. 
This suggests to consider an auxiliary problem by removing the problematic constraint 
$G(x)*H(y)=0$ and including other ones that ensure the null product. We then define 
the $I$-{\em Tightened} Nonlinear Problem at $(\bar{x},\bar{y})$ by
\begin{equation}
\label{aw_prob:tight}
\begin{array}{cl}
\displaystyle\mathop{\rm minimize }_{x,y} & f(x)  \\
{\rm subject\ to } & g(x)\leq 0, h(x)=0, \\
& \theta(y)\leq 0, \\
& G_i(x)=0,\; i \in I, \\
& H_i(y)\leq 0,\; i\in I_{0+}(\bar{x},\bar{y})\cup I_{01}(\bar{x},\bar{y}), \\
& H_i(y)=0,\; i\in I_{00}(\bar{x},\bar{y})\cup I_{\pm0}(\bar{x},\bar{y}), \\
& \tilde{H}(y)\leq 0.
\end{array}
\end{equation}
This problem will be also indicated by TNLP$_{I}(\bar{x},\bar{y})$ 
and, when there is no chance for ambiguity, it will be referred simply to 
as {\em tightened problem}. Note that we tighten only 
those constraints that are involved in the complementarity constraint $G(x)*H(y)=0$, by 
incorporating the equality constraints $G_i$'s and converting the active inequalities 
$H_i$'s into equalities.

The following lemma is a straightforward consequence of the definition of the tightened 
problem TNLP$_{I}(\bar{x},\bar{y})$. 
\begin{lemma}
\label{aw_lm:tnlp1}
Consider the tightened problem (\ref{aw_prob:tight}). Then, 
\begin{enumerate}
\item\label{aw_inactiveH} the inequalities defined by $H_i$, 
$i\in I_{0+}(\bar{x},\bar{y})\cup I_{01}(\bar{x},\bar{y})$, are inactive at $(\bar{x},\bar{y})$;
\item\label{aw_xybar_tnlp} $(\bar{x},\bar{y})$ is feasible for TLNP$_{I}(\bar{x},\bar{y})$; 
\item every feasible point of (\ref{aw_prob:tight}) is feasible for (\ref{aw_prob:relax1});
\item if $(\bar{x},\bar{y})$ is a global (local) minimizer of (\ref{aw_prob:relax1}), 
then it is also a global (local) minimizer of TNLP$_{I}(\bar{x},\bar{y})$.
\end{enumerate}
\end{lemma}

The Lagrangian function associated with TNLP$_{I}(\bar{x},\bar{y})$ is the function 
$$
\mathcal{L}_{I}:\R^n\times\R^n\times\R^m\times\R^p\times\R\times\R^{|I|}
\times\R^n\times\R^n\to\R
$$ 
given by 
\begin{align*}
\mathcal{L}_{I}(x,y,\lambda^g,\lambda^h,\lambda^{\theta},\lambda_I^G,\lambda^H,
\lambda^{\tilde{H}}) = & \, L(x,\lambda^g,\lambda^h)+
\lambda^{\theta}\theta(y)+(\lambda_I^G)^TG_{I}(x) \\ 
& +(\lambda^H)^TH(y)+(\lambda^{\tilde{H}})^T\tilde{H}(y),
\end{align*}
where $L$ is the Lagrangian defined in (\ref{aw_lagr:nlp}).

Note that the tightened problem, and hence its Lagrangian, depends on the index set $I$, which 
in turn depends on the point $(\bar{x},\bar{y})$. It should be also noted that 
\begin{equation}
\label{aw_nablaLIx}
\nabla_{x,y}\mathcal{L}_{I}(x,y,\lambda)=
\left(\begin{array}{c} \nabla_{x}L(x,\lambda^g,\lambda^h)+\sum_{i\in I}\lambda_i^Ge_i 
\vspace{3pt} \\ -\lambda^{\theta}e-\lambda^H+\lambda^{\tilde{H}} \end{array}\right).
\end{equation}

The weaker stationarity concept proposed in \cite{KrulikovskiRibeiroSachine20aX} 
is presented below. 
\begin{definition}
\label{aw_def:wstat_xy}
Consider a feasible point $(\bar{x},\bar{y})$ of the relaxed problem (\ref{aw_prob:relax1}) 
and a set of indices $I$ satisfying (\ref{aw_indexI}). We say that $(\bar{x},\bar{y})$ 
is $I$-weakly stationary ($W_{I}$-stationary) for this problem if there exists a vector 
$$
\lambda=(\lambda^g,\lambda^h,\lambda^{\theta},\lambda_I^G,\lambda^H,\lambda^{\tilde{H}})\in
\R_+^m\times\R^p\times\R_+\times\R^{|I|}\times\R^n\times\R_+^n
$$ 
such that
\begin{enumerate}
\item\label{aw_lagrangian0} $\nabla_{x,y}{\cal{L}}_{I}(\bar{x},\bar{y},\lambda)=0$;
\item\label{aw_g0} $(\lambda^g)^Tg(\bar{x})=0$;
\item\label{aw_theta0} $\lambda^{\theta}\theta(\bar{y})=0$;
\item\label{aw_H0} $\lambda^H_i=0$ for all $i\in I_{0+}(\bar{x},\bar{y})
\cup I_{01}(\bar{x},\bar{y})$;
\item\label{aw_htilde0} $(\lambda^{\tilde{H}})^T\tilde{H}(\bar{y})=0$.
\end{enumerate}    
\end{definition}

\begin{remark}
\label{aw_rm:wstat}
The first item of Definition \ref{aw_def:wstat_xy} is nothing else than the gradient 
KKT condition for the tightened problem (\ref{aw_prob:tight}). Items (\ref{aw_g0}), (\ref{aw_theta0}) 
and (\ref{aw_htilde0}) represent the standard KKT complementarity conditions for the 
inequality constraints $g(x)\leq 0$, $\theta(y)\leq 0$ and $\tilde{H}(y)\leq 0$, 
respectively, of the tightened problem. Item (\ref{aw_H0}) also represents KKT 
complementarity conditions for the constraints $H_i(y)\leq 0$, 
$i\in I_{0+}(\bar{x},\bar{y})\cup I_{01}(\bar{x},\bar{y})$, in view of 
Lemma \ref{aw_lm:tnlp1}(\ref{aw_inactiveH}).
\end{remark}

As an immediate consequence of Remark \ref{aw_rm:wstat} we have the following characterization 
of $W_{I}$-stationarity for the relaxed problem in terms of stationarity for the 
tightened problem. 
\begin{proposition}
\label{aw_prop:wstat_kkttnlp}
Let $(\bar{x},\bar{y})$ be a feasible point of the relaxed problem (\ref{aw_prob:relax1}). Then, 
$(\bar{x},\bar{y})$ is $W_{I}$-stationary if and only if it is a KKT point for the tightened 
problem (\ref{aw_prob:tight}).
\end{proposition}

Note that in view of Proposition \ref{aw_prop:wstat_kkttnlp} we could have defined 
$W_{I}$-stationarity simply as KKT for the tightened problem (\ref{aw_prob:tight}). Nevertheless, 
we prefer as in Definition \ref{aw_def:wstat_xy} in order to have its last condition (\ref{aw_H0}) 
explicitly, instead of hidden in the complementarity condition. This way of stating weak 
stationarity is also similar to that used in the MPCC setting, see 
\cite{AndreaniHaeserSecchinSilva,FlegelKanzow}.

In the next result we justify why Definition \ref{aw_def:wstat_xy} is considered a weaker 
stationarity concept for the relaxed problem. 
\begin{theorem}
\label{aw_th:kkt_wstat}
Suppose that $(\bar{x},\bar{y})$ is a KKT point for the relaxed problem (\ref{aw_prob:relax1}). 
Then $(\bar{x},\bar{y})$ is $W_{I}$-stationary.
\end{theorem}

At this point we could ask if $W_{I}$-stationarity, being weaker than KKT, is a 
necessary optimality condition. That is, can we ensure that a minimizer of the relaxed 
problem is $W_{I}$-stationary for some index set $I$ satisfying (\ref{aw_indexI})? 
The answer is no, as illustrated in the following example.

\begin{example}
\label{aw_ex:minnotWI}
Consider the MPCaC and the associated relaxed problem given below.
$$
\begin{array}{lr}
\begin{array}{cl}
\displaystyle\mathop{\rm minimize }_{x\in\R^3} & x_1  \\
{\rm subject\ to } & (1-x_1)^3+x_3^2\leq 0, \\
& \|x\|_0\leq 2, \\ 
& \\ 
&
\end{array}
\hspace{.5cm}
&
\begin{array}{cl}
\displaystyle\mathop{\rm minimize }_{x,y\in\R^3} & x_1 \\
{\rm subject\ to } & (1-x_1)^3+x_3^2\leq 0, \\
& y_1+y_2+y_3\geq 1, \\
& x_iy_i=0,\; i=1,2,3, \\
& 0\leq y_i \leq 1, \; i=1,2,3.
\end{array}
\end{array}
$$
The point $x^*=(1,0,0)$ is a global solution of MPCaC and $(x^*,y^*)$, with  
$y^*=(0,1,0)$, is a global solution of the relaxed problem.
For the points $x^*$ and $(x^*,y^*)$ we have 
$$
I_{0}=\{2,3\},\;\; I_{01}=\{2\},\;\; I_{\pm 0}=\{1\},\;\; I_{00}=\{3\}
\quad\mbox{and}\quad I_{0+}=\emptyset.
$$
So, there are two choices for $I$ that satisfy (\ref{aw_indexI}): $I'=\{2\}$ or $I''=\{2,3\}$. 
Let us analyze each one of them. 

For $I=I'$ we have 
\begin{align*}
\nabla_{x}L(x^*,\lambda^g)+\sum_{i\in I}\lambda_i^Ge_i=
\left(\begin{array}{c} 1 \\ 0 \\ 0 \end{array}\right)+
\left(\begin{array}{c} 0 \\ \lambda_2^G \\ 0 \end{array}\right).
\end{align*}
Since this expression does not vanish, taking into account (\ref{aw_nablaLIx}), we see that 
the pair $(x^*,y^*)$ is not $W_{I}$-stationary. 

Now, for $I=I''$ we have 
\begin{align*}
\nabla_{x}L(x^*,\lambda^g)+\sum_{i\in I}\lambda_i^Ge_i=
\left(\begin{array}{c} 1 \\ 0 \\ 0 \end{array}\right)+
\left(\begin{array}{c} 0 \\ \lambda_2^G \\ \lambda_3^G \end{array}\right).
\end{align*}
Again, the above expression does not vanish and then $(x^*,y^*)$ is not $W_{I}$-stationary. 
\end{example}

In view of Example \ref{aw_ex:minnotWI} and motivated to find a necessary optimality 
condition for MPCaC problems, we propose in the next section the concept of 
approximate weak stationarity, which will be satisfied at every minimizer, 
independently of any constraint qualification. 

\section{Sequential optimality conditions for MPCaC}
\label{aw_sec:aw}
In order to define our sequential optimality condition, consider the function 
$$
\mathcal{L}:\R^n\times\R^n\times\R^m\times\R^p\times\R\times\R^n
\times\R^n\times\R^n\to\R
$$ 
defined by  
\begin{align*}
\mathcal{L}(x,y,\lambda^g,\lambda^h,\lambda^{\theta},\lambda^G,\lambda^H,
\lambda^{\tilde{H}}) = & \, L(x,\lambda^g,\lambda^h)+
\lambda^{\theta}\theta(y)+(\lambda^G)^TG(x) \\ 
& +(\lambda^H)^TH(y)+(\lambda^{\tilde{H}})^T\tilde{H}(y),
\end{align*}
where $L$ is the Lagrangian defined in (\ref{aw_lagr:nlp}).

Note that ${\cal L}$ resembles the Lagrangian ${\cal L}_I$, associated with 
TNLP$_I(\bar{x},\bar{y})$. The only difference is that the term 
$(\lambda_I^G)^TG_{I}(x)$ was replaced by $(\lambda^G)^TG(x)$. Here it will be 
convenient to see that 
\begin{align}
\nabla_{x,y}\mathcal{L}(x,y,\lambda)=
\left(\begin{array}{c} \nabla_{x}L(x,\lambda^g,\lambda^h)+
\displaystyle\sum_{i=1}^{n}\lambda_i^G\nabla G_i(x) \\ \lambda^{\theta}\nabla\theta(y) 
+\displaystyle\sum_{i=1}^{n}\lambda_i^H \nabla H_i(y)
+\sum_{i=1}^{n}\lambda_i^{\tilde{H}}\nabla \tilde{H}_i(y) 
\end{array}\right). \label{aw_nablaL0}
\end{align}

\begin{definition}
\label{aw_def:aw}
Let $(\bar{x},\bar{y})$ be a feasible point of the relaxed problem (\ref{aw_prob:relax1}). 
We say that $(\bar{x},\bar{y})$ is Approximately Weakly stationary ($AW$-stationary) 
for this problem if there exist sequences $(x^k,y^k)\subset\R^n\times\R^n$ and 
$$
(\lambda^k)=\big(\lambda^{g,k},\lambda^{h,k},\lambda^{\theta,k},\lambda^{G,k},
\lambda^{H,k},\lambda^{\tilde{H},k}\big)\subset\R_+^m\times\R^p\times\R_+
\times\R^n\times\R^n\times\R_+^n
$$
such that
\begin{enumerate}
\item\label{aw_xkyk} $(x^k,y^k)\to(\bar{x},\bar{y})$; \vspace{2pt}
\item\label{aw_gradLtnlp} $\nabla_{x,y}{\cal{L}}(x^k,y^k,\lambda^k)\to 0$; \vspace{2pt}
\item\label{aw_glbdg} $\min\{-g(x^k),\lambda^{g,k}\}\to 0$; \vspace{2pt}
\item\label{aw_thlbdth} $\min\{-\theta(y^k),\lambda^{\theta,k}\}\to 0$; 
\vspace{2pt}
\item\label{aw_GlbdG} $\min\{|G_i(x^k)|,|\lambda_i^{G,k}|\}\to 0$ for all 
$i=1,\ldots,n$;  \vspace{2pt}
\item\label{aw_HlbdH} $\min\{-H_i(y^k),|\lambda_i^{H,k}|\}\to 0$ for all 
$i=1,\ldots,n$; \vspace{2pt}
\item\label{aw_HtillbdHtil} $\min\{-\tilde{H}(y^k),\lambda^{\tilde{H},k}\}\to 0$.
\end{enumerate} 
\end{definition}

\begin{remark}
\label{aw_rm:awstat1}
Definition~\ref{aw_def:aw} resembles AKKT condition, where (\ref{aw_glbdg}), (\ref{aw_thlbdth}) 
and (\ref{aw_HtillbdHtil}) represent the approximate complementarity conditions for 
the inequality constraints $g(x)\leq 0$, $\theta(y)\leq 0$ and $\tilde{H}(y)\leq 0$, 
respectively and (\ref{aw_HlbdH}) is related to the last complementarity condition in 
W$_I$-stationarity. As a matter of fact, $AW$-stationarity is equivalent to AKKT for 
TNLP$_{I_0}$, as we shall see ahead in Theorem \ref{aw_th:aw_akkt}.
\end{remark}

Let us review Example \ref{aw_ex:minnotWI} in light of the above definition. We have seen 
that the minimizer is not $W_{I}$-stationary, but now we can see that it is 
$AW$-stationary. 

\begin{example}
\label{aw_ex:minnotWI_AW}
Consider the problem given in Example \ref{aw_ex:minnotWI}.
We claim that the global solution of the relaxed problem, $(x^*,y^*)$, is 
$AW$-stationary. Indeed, 
consider the sequences $(x^k,y^k)\subset\R^3\times\R^3$ and 
$$
(\lambda^k)=\big(\lambda^{g,k},\lambda^{\theta,k},\lambda^{G,k},
\lambda^{H,k},\lambda^{\tilde{H},k}\big)\subset\R_+^3\times\R_+\times\R^3\times\R^3\times\R_+^3
$$
defined by $x^k=(1+1/k,0,0)$, $y^k=(0,1,0)$, $\lambda^{g,k}=k^2/3$, $\lambda^{\theta,k}=0$ and 
$\lambda^{G,k}=\lambda^{H,k}=\lambda^{\tilde{H},k}=0$. Then, we have $(x^k,y^k)\to(x^*,y^*)$ 
and 
\begin{align*}
\nabla_{x}L(x^k,\lambda^{g,k})+\sum_{i=1}^{n}\lambda_i^{G,k}\nabla G_i(x^k)=
\left(\begin{array}{c} 1-3\lambda^{g,k}(1-x_1^k)^2 \\ 0 \\ 2\lambda^{g,k}x_3^k 
\end{array}\right)=0.
\end{align*}
So, in view of (\ref{aw_nablaL0}), we obtain the first two items of Definition \ref{aw_def:aw}. 
Now, note that 
$g(x^k)\to g(x^*)=0$ and $\theta(y^k)\to\theta(y^*)=0$, 
which in turn imply that 
$$
\min\{-g(x^k),\lambda^{g,k}\}\to 0\quad\mbox{and}\quad\min\{-\theta(y^k),\lambda^{\theta,k}\}\to 0,
$$
giving items (\ref{aw_glbdg}) and 
(\ref{aw_thlbdth}). The relation $\min\{|G_i(x^k)|,|\lambda_i^{G,k}|\}\to 0$ is immediate. 
Besides, since $\tilde{H}(y^k)\to\tilde{H}(y^*)\leq 0$, $\lambda^{\tilde{H},k}=0$, 
$H(y^k)\to H(y^*)\leq 0$ and $\lambda^{H,k}=0$, 
we have $\min\{-\tilde{H}(y^k),\lambda^{\tilde{H},k}\}\to 0$ and 
$\min\{-H_i(y^k),|\lambda_i^{H,k}|\}\to 0$, obtaining items 
(\ref{aw_GlbdG}), (\ref{aw_HlbdH}) and  (\ref{aw_HtillbdHtil}).
\end{example}

Now we shall prove that the above example reflects a general result, that is, 
every minimizer of an MPCaC problem is $AW$-stationary. 
We start the theoretical analysis with two simple facts. 
The first one says that the expression $\sum_{i=1}^{n}\lambda_i^{G,k}\nabla G_i(x^k)$ 
could be replaced by 
$\sum_{i\in I_0}\lambda_i^{G,k}\nabla G_i(x^k)$. The second fact states 
that $AW$-stationarity is weaker than $W_{I}$-stationarity, and consequently weaker 
than KKT, in view of Theorem~\ref{aw_th:kkt_wstat}.

\begin{lemma}
\label{aw_lm:awI0}
Let $(\bar{x},\bar{y})$ be an $AW$-stationary point for the relaxed 
problem (\ref{aw_prob:relax1}), with corresponding sequences $(x^k,y^k)$ and 
$(\lambda^k)$. Then,  
$$
\nabla_xL(x^k,\lambda^{g,k},\lambda^{h,k})+
\sum_{i\in I_0}\lambda_i^{G,k}\nabla G_i(x^k)\to 0.
$$
\end{lemma}
\beginproof
In view of (\ref{aw_nablaL0}), we have, in particular, 
\begin{equation}
\label{aw_eq_awI01}
\nabla_xL(x^k,\lambda^{g,k},\lambda^{h,k})+
\sum_{i=1}^{n}\lambda_i^{G,k}\nabla G_i(x^k)\to 0.
\end{equation}
For $i\notin I_0$, we have $\lim_{k\to\infty}G_i(x^k)=G_i(\bar{x})=\bar{x}_i\neq 0$. 
Therefore, we can assume without loss of generality that there exists $\epsilon>0$ 
such that $|G_i(x^k)|\geq\epsilon$ for all $k$. Since 
$\min\{|G_i(x^k)|,|\lambda_i^{G,k}|\}\to 0$, we obtain  
$|\lambda_i^{G,k}|\to 0$ and hence, 
$$
\sum_{i\notin I_0}\lambda_i^{G,k}\nabla G_i(x^k)\to 0.
$$
By subtracting this from (\ref{aw_eq_awI01}), we conclude the proof.
\endproof

\begin{lemma}
\label{aw_lm:ws_aw}
Let $(\bar{x},\bar{y})$ be a $W_{I}$-stationary point for the relaxed 
problem (\ref{aw_prob:relax1}), in the sense of Definition \ref{aw_def:wstat_xy}. Then 
$(\bar{x},\bar{y})$ is $AW$-stationary for this problem.
\end{lemma}
\beginproof
Consider a vector 
$$
\lambda=(\lambda^g,\lambda^h,\lambda^{\theta},\lambda_I^G,\lambda^H,
\lambda^{\tilde{H}})\in\R_+^m\times\R^p\times\R_+\times\R^{|I|}\times\R^n\times\R_+^n
$$ 
satisfying Definition \ref{aw_def:wstat_xy}. Then, the (constant) sequences 
$(x^k,y^k)\subset\R^n\times\R^n$ and 
$$
(\lambda^k)=\big(\lambda^{g,k},\lambda^{h,k},\lambda^{\theta,k},\lambda^{G,k},
\lambda^{H,k},\lambda^{\tilde{H},k}\big)\subset\R_+^m\times\R^p\times\R_+
\times\R^n\times\R^n\times\R_+^n,
$$
defined by 
$$
(x^k,y^k)=(\bar{x},\bar{y})\,,\ \big(\lambda^{g,k},\lambda^{h,k},\lambda^{\theta,k},
\lambda_I^{G,k},\lambda^{H,k},\lambda^{\tilde{H},k}\big)=(\lambda^g,\lambda^h,
\lambda^{\theta},\lambda_I^G,\lambda^H,\lambda^{\tilde{H}})
$$
and $\lambda_i^{G,k}=0$ for $i\notin I$ and $k\in\mathbb{N}$, satisfy 
Definition \ref{aw_def:aw}. 
\endproof

\begin{remark}
\label{aw_rm:awstat2}
We point out here that, in contrast to W$_I$-stationarity, which conveniently depends 
on the set $I$, our sequential optimality condition is independent of any set $I$. 
This is a desirable feature since $AW$-stationarity has a certain amount of algorithmic 
appeal. In practice, one is able to use such conditions as a stopping criterion for 
an algorithm designed to solve the MPCaC problem. 
\end{remark}

Before proving our main sequential optimality results, let us see some preliminary 
lemmas. To this end, consider the augmented problem 
\begin{equation}
\label{aw_prob:augm}
\begin{array}{cl}
\displaystyle\mathop{\rm minimize }_{x,y,w} & f(x)  \\
{\rm subject\ to } & g(x)\leq 0, h(x)=0, \\
& \theta(y)\leq 0, \\
& w^G-G(x)=0,\, w^H+H(y)=0, \\
& \tilde{H}(y)\leq 0, \\
& w\in W,
\end{array}
\end{equation}
where $W=\{w=\big(w^G,w^H\big)\in\R^n\times\R_+^n\mid w^G*w^H=0\}$. 

This problem will be crucial in the analysis. In the next two lemmas we establish 
the equivalence between the relaxed problem (\ref{aw_prob:relax1}) and this augmented 
problem. Moreover, there is a suitable reason to write the constraints 
$H(y)\leq 0$ and $G(x)*H(y)=0$ of (\ref{aw_prob:relax1}) in the format 
$w\in W$. Such a strategy will enable us to apply Lemma \ref{aw_lm:compl1} 
to obtain Guignard constraint qualification for an auxiliary problem ahead.

\begin{lemma}
\label{aw_lm:sol_relax_sol_augm}
Let $(x^*,y^*)$ be a local (global) minimizer of the relaxed problem (\ref{aw_prob:relax1}). 
Given $w^*\in\R^n\times\R^n$, if the point $(x^*,y^*,w^*)$ is 
feasible for the augmented problem~(\ref{aw_prob:augm}), then it is a local 
(global) minimizer of this problem. In particular, this holds for 
$w^*=\big(G(x^*),-H(y^*)\big)$.
\end{lemma}
\beginproof
First, the relation between local minimizers. In view of the equivalence of norms, 
we consider $\|\cdot\|_{\infty}$, for convenience. 
By hypothesis, there exists $\delta>0$ such that if $(x,y)$ is feasible 
for (\ref{aw_prob:relax1}) and $\|(x,y)-(x^*,y^*)\|_{\infty}\leq\delta$, 
then $f(x^*)\leq f(x)$. Suppose that $(x^*,y^*,w^*)$ is feasible for the 
problem (\ref{aw_prob:augm}) and consider an arbitrary feasible point $(x,y,w)$ for 
this problem such that $\|(x,y,w)-(x^*,y^*,w^*)\|_{\infty}\leq\delta$. 
Then, the pair $(x,y)$ is feasible for (\ref{aw_prob:relax1}) and 
$\|(x,y)-(x^*,y^*)\|_{\infty}\leq\delta$. Hence, $f(x^*)\leq f(x)$ and, 
therefore, $(x^*,y^*,w^*)$ is a local minimizer of (\ref{aw_prob:augm}). Note 
that $(x^*,y^*,w^*)$, with $w^*=\big(G(x^*),-H(y^*)\big)$, is trivially feasible. 
Finally, if we ignore the neighborhoods in the argument above, we obtain the 
relation between global minimizers.
\endproof

For the sake of completeness we prove below the converse of 
Lemma \ref{aw_lm:sol_relax_sol_augm}.

\begin{lemma}
\label{aw_lm:sol_augm_sol_relax}
Let $(x^*,y^*,w^*)$ be a local (global) minimizer of (\ref{aw_prob:augm}). 
Then $(x^*,y^*)$ is a local (global) minimizer of (\ref{aw_prob:relax1}). 
\end{lemma}
\beginproof
By the feasibility of $(x^*,y^*,w^*)$ we have that $(x^*,y^*)$ is feasible for 
(\ref{aw_prob:relax1}), 
\begin{equation}
\label{aw_eq:wGH1}
(w^*)^G=G(x^*) \quad\mbox{and}\quad (w^*)^H=-H(y^*).
\end{equation}
Consider $\delta_1>0$ such that $f(x^*)\leq f(x)$ for all feasible point $(x,y,w)$ of 
(\ref{aw_prob:augm}), satisfying $\|(x,y,w)-(x^*,y^*,w^*)\|_{\infty}\leq\delta_1$. 
Let $\delta_2>0$ be such that 
\begin{equation}
\label{aw_eq:wGH2}
\|G(x)-G(x^*)\|_{\infty}\leq\delta_1 \quad\mbox{and}\quad 
\|H(y)-H(y^*)\|_{\infty}\leq\delta_1
\end{equation}
for all $(x,y)\in\R^n\times\R^n$ with $\|(x,y)-(x^*,y^*)\|_{\infty}\leq\delta_2$. Define 
$\delta=\min\{\delta_1,\delta_2\}$ and take $(x,y)$, feasible for (\ref{aw_prob:relax1}),  
such that $\|(x,y)-(x^*,y^*)\|_{\infty}\leq\delta$. Thus we have (\ref{aw_eq:wGH2}), which in 
view of (\ref{aw_eq:wGH1}) can be rewritten as $\|w-w^*\|_{\infty}\leq\delta_1$, with 
$w=\big(G(x),-H(y)\big)$. Therefore, $(x,y,w)$ is feasible for (\ref{aw_prob:augm}) 
and 
$$
\|(x,y,w)-(x^*,y^*,w^*)\|_{\infty}\leq\delta_1,
$$ 
implying that $f(x^*)\leq f(x)$. 

Now, let us see the global optimality. So, assume that $(x^*,y^*,w^*)$ is a global minimizer 
of (\ref{aw_prob:augm}). Then $(x^*,y^*)$ is feasible for (\ref{aw_prob:relax1}). 
Furthermore, given an arbitrary feasible point $(x,y)$, we have that 
$(x,y,w)$, with $w=\big(G(x),-H(y)\big)$, is feasible for (\ref{aw_prob:augm}). 
Therefore,  $f(x^*)\leq f(x)$.
\endproof

\begin{lemma}
\label{aw_lm:sol_relax_sol_prox}
Suppose that $(x^*,y^*)$ is a local minimizer of the relaxed problem (\ref{aw_prob:relax1}). 
Then, given an arbitrary norm $\|\cdot\|$, there exists $\delta>0$ such that 
$(x^*,y^*,w^*)$, with $w^*=\big(G(x^*),-H(y^*)\big)$, is the unique global minimizer 
of the problem 
\begin{equation}
\label{aw_prob:prox}
\begin{array}{cl}
\displaystyle\mathop{\rm minimize }_{x,y,w} & f(x)+\dfrac{1}{2}\|(x,y)-(x^*,y^*)\|_2^2  \\
{\rm subject\ to } & g(x)\leq 0, h(x)=0, \\
& \theta(y)\leq 0, \\
& w^G-G(x)=0,\, w^H+H(y)=0, \\
& \tilde{H}(y)\leq 0, \\
& w\in W, \\
& \|(x,y,w)-(x^*,y^*,w^*)\|\leq\delta.
\end{array}
\end{equation}
\end{lemma}
\beginproof
By Lemma \ref{aw_lm:sol_relax_sol_augm}, we have that $(x^*,y^*,w^*)$ is a local minimizer 
of (\ref{aw_prob:augm}). 
Consider $\delta>0$ such that if $(x,y,w)$ is feasible for (\ref{aw_prob:augm}) and 
\begin{equation}
\label{aw_eq:xywlocal}
\|(x,y,w)-(x^*,y^*,w^*)\|\leq\delta,
\end{equation}
then $f(x^*)\leq f(x)$. 
Note that $(x^*,y^*,w^*)$ is feasible for (\ref{aw_prob:prox}). Moreover, given any 
feasible point $(x,y,w)$, of (\ref{aw_prob:prox}), we have that it is also feasible for 
(\ref{aw_prob:augm}) and satisfies (\ref{aw_eq:xywlocal}). Hence, 
$$
f(x^*)+\dfrac{1}{2}\|(x^*,y^*)-(x^*,y^*)\|_2^2=f(x^*)\leq f(x)\leq f(x)+
\dfrac{1}{2}\|(x,y)-(x^*,y^*)\|_2^2,
$$
proving that $(x^*,y^*,w^*)$ is a global minimizer of (\ref{aw_prob:prox}). 

Now, suppose that $(\bar{x},\bar{y},\bar{w})$ is also a global minimizer 
of (\ref{aw_prob:prox}). Then, 
$$
f(\bar{x})+\dfrac{1}{2}\|(\bar{x},\bar{y})-(x^*,y^*)\|_2^2\leq f(x^*)+
\dfrac{1}{2}\|(x^*,y^*)-(x^*,y^*)\|_2^2=f(x^*)\leq f(\bar{x}),
$$
where the last inequality follows from the fact that $(\bar{x},\bar{y},\bar{w})$ is 
feasible for (\ref{aw_prob:augm}) and satisfies (\ref{aw_eq:xywlocal}). Therefore, 
$(\bar{x},\bar{y})=(x^*,y^*)$, and hence 
$$
\bar{w}=\big(G(\bar{x}),-H(\bar{y})\big)=
\big(G(x^*),-H(y^*)\big)=w^*,
$$
proving the uniqueness.
\endproof

The next result shows that our stationarity concept, given in Definition \ref{aw_def:aw}, 
is a legitimate optimality condition, independently of any constraint qualification. 
This is a requirement for them to be useful in the analysis of algorithms.

\begin{theorem}
\label{aw_th:aw}
If $(x^*,y^*)$ is a local minimizer of the relaxed problem (\ref{aw_prob:relax1}), then 
it is an $AW$-stationary point, in the sense of Definition \ref{aw_def:aw}. 
\end{theorem}
\beginproof
Defining $w^*=\big(G(x^*),-H(y^*)\big)$, we conclude from 
Lemma \ref{aw_lm:sol_relax_sol_prox} that there exists $\delta>0$ such that the point 
$(x^*,y^*,w^*)$ is the unique global minimizer of the problem~(\ref{aw_prob:prox}), 
with $\|\cdot\|_2$ in the last constraint. Define the (partial) infeasibility 
measure associated with this problem as 
\begin{equation*}
\begin{array}{rcl}
\varphi(x,y,w) & = & \dfrac{1}{2}\Big(\|g^+(x)\|_2^2+\|h(x)\|_2^2+\|\theta^+(y)\|_2^2+
\|w^G-G(x)\|_2^2 \vspace{3pt} \\ 
& & +\|w^H+H(y)\|_2^2+\|\tilde H^+(y)\|_2^2\Big),
\end{array}
\end{equation*}
consider a sequence $\rho_k\to\infty$ and let $(x^k,y^k,w^k)$ be a global minimizer 
of the penalized problem 
\begin{equation}
\label{aw_prob:penalized}
\begin{array}{cl}
\displaystyle\mathop{\rm minimize }_{x,y,w} & f(x)+\dfrac{1}{2}\|(x,y)-(x^*,y^*)\|_2^2 +
\rho_k\varphi(x,y,w) \\
{\rm subject\ to } & w\in W, \\
& \|(x,y,w)-(x^*,y^*,w^*)\|_2^2\leq\delta^2,
\end{array}
\end{equation}
which is well defined because the objective function is continuous and the feasible 
set is compact. Since $\|(x^k,y^k,w^k)-(x^*,y^*,w^*)\|_2\leq\delta$, we can assume 
without loss of generality that the sequence $(x^k,y^k,w^k)$ converges to some point 
$(\bar{x},\bar{y},\bar{w})$. We claim that $(\bar{x},\bar{y},\bar{w})=(x^*,y^*,w^*)$. 
Note first that $(x^*,y^*,w^*)$ is feasible for (\ref{aw_prob:penalized}) and 
$\varphi(x^*,y^*,w^*)=0$. So, by the optimality of $(x^k,y^k,w^k)$ we have 
\begin{equation}
\label{aw_eq:penaliz}
f(x^k)+\dfrac{1}{2}\|(x^k,y^k)-(x^*,y^*)\|_2^2+
\rho_k\varphi(x^k,y^k,w^k)\leq f(x^*),
\end{equation}
implying that $\varphi(x^k,y^k,w^k)\to 0$, because $\rho_k\to\infty$. 
This in turn implies that $\varphi(\bar{x},\bar{y},\bar{w})=0$, giving 
$g^+(\bar{x})=0$, $h(\bar{x})=0$, $\theta^+(\bar{y})=0$, 
$\bar{w}^G=G(\bar{x})$, $\bar{w}^H=-H(\bar{y})$ and $\tilde{H}^+(\bar{y})=0$. 
Moreover, as the sequence $(x^k,y^k,w^k)$ is feasible for (\ref{aw_prob:penalized}), 
its limit point $(\bar{x},\bar{y},\bar{w})$ satisfies $\bar{w}\in W$, because $W$ 
is a closed set, and $\|(\bar{x},\bar{y},\bar{w})-(x^*,y^*,w^*)\|\leq\delta$. 
Therefore, $(\bar{x},\bar{y},\bar{w})$ is feasible for (\ref{aw_prob:prox}). 
Furthermore, from (\ref{aw_eq:penaliz}) we obtain
$$
f(x^k)+\dfrac{1}{2}\|(x^k,y^k)-(x^*,y^*)\|_2^2\leq f(x^*).
$$
Taking the limit it follows that 
$$
f(\bar{x})+\dfrac{1}{2}\|(\bar{x},\bar{y})-(x^*,y^*)\|_2^2\leq f(x^*),
$$
which means that $(\bar{x},\bar{y},\bar{w})$ is optimal for (\ref{aw_prob:prox}). 
By the uniqueness of the optimal solution of this problem, we conclude that 
$(\bar{x},\bar{y},\bar{w})=(x^*,y^*,w^*)$, proving the claim. As consequence, 
we have the first item of Definition \ref{aw_def:aw}. 

In order to prove the next item, let us see first that a constraint 
qualification holds at the minimizer $(x^k,y^k,w^k)$. Since 
$(x^k,y^k,w^k)\to(x^*,y^*,w^*)$, we may assume without loss of generality that 
$\|(x^k,y^k,w^k)-(x^*,y^*,w^*)\|_2<\delta$ for all $k$. That is, the inequality 
constraint in the problem (\ref{aw_prob:penalized}) is inactive at the minimizer. 
By Lemma \ref{aw_lm:inactive}, the tangent and linearized cones at 
this point are the ones taking into account only the constraints in $w\in W$, namely, 
\begin{equation}
\label{aw_eq:wW}
-w^H\leq 0\quad\mbox{and}\quad w^G*w^H=0.
\end{equation}
Thus, in view of Lemmas \ref{aw_lm:withoutx} and \ref{aw_lm:compl1}, Guignard constraint 
qualification holds at $(x^k,y^k,w^k)$. 
This implies that it satisfies the KKT conditions, which means that there exist 
multipliers $\mu^{H,k}\in\R_+^n$ and $\mu^{0,k}\in\R^n$, associated with the 
constraints in $w\in W$, such that 
\begin{subequations}
\begin{align}
\nabla f(x^k)+(x^k-x^*)+\rho_k\nabla_{x}\varphi(x^k,y^k,w^k)=0 \label{aw_gradLpx} \\ 
(y^k-y^*)+\rho_k\nabla_{y}\varphi(x^k,y^k,w^k)=0 \label{aw_gradLpy} \\ 
\rho_k\nabla_{w^G}\varphi(x^k,y^k,w^k)+\mu^{0,k}*w^{H,k}=0 \label{aw_gradLpwG}  \\ 
\rho_k\nabla_{w^H}\varphi(x^k,y^k,w^k)-\mu^{H,k}+\mu^{0,k}*w^{G,k}=0 \label{aw_gradLpwH}  \\ 
\mu^{H,k}*w^{H,k}=0. \label{aw_complw}
\end{align}
\end{subequations}
Noting that the partial gradients of $\varphi$ are given by 
\begin{subequations}
\begin{align}
\nabla_{x}\varphi(x,y,w)=\nabla g(x)g^+(x)+\nabla h(x)h(x)+
\nabla G(x)\big(G(x)-w^G\big), \label{aw_gradphix} \\ 
\nabla_{y}\varphi(x,y,w)=\theta^+(y)\nabla\theta(y)+\nabla\tilde H(y)\tilde H^+(y)+
\nabla H(y)\big(w^H+H(y)\big), \label{aw_gradphiy} \\ 
\nabla_{w^G}\varphi(x,y,w)=w^G-G(x)\quad\mbox{and}\quad 
\nabla_{w^H}\varphi(x,y,w)=w^H+H(y) \label{aw_gradphiw}
\end{align}
\end{subequations}
and defining $\lambda^k$ as  
$$
\begin{array}{c}
\lambda^{g,k}=\rho_kg^+(x^k),\ \lambda^{h,k}=\rho_kh(x^k),\ \lambda^{\theta,k}=
\rho_k\theta^+(y^k), 
\vspace{6pt} \\ {\lambda}^{G,k}=\rho_k\big(G(x^k)-w^{G,k}\big),
\ \lambda^{H,k}=\rho_k\big(w^{H,k}+H(y^k)\big),
\ \lambda^{\tilde{H},k}=\rho_k\tilde{H}^+(y^k),
\end{array}
$$
we see immediately that $\lambda^{g,k}\geq 0$, $\lambda^{\theta,k}\geq 0$ and 
$\lambda^{\tilde{H},k}\geq 0$. Moreover, using (\ref{aw_gradLpx}) and (\ref{aw_gradphix}), 
we obtain 
$$
\nabla_{x}{\cal L}(x^k,y^k,\lambda^k)=
\nabla f(x^k)+\rho_k\nabla_{x}\varphi(x^k,y^k,w^k)=x^*-x^k\to 0.
$$
Furthermore, from (\ref{aw_gradLpy}) and (\ref{aw_gradphiy}), we have 
$$
\nabla_{y}{\cal L}(x^k,y^k,\lambda^k)=
\rho_k\nabla_{y}\varphi(x^k,y^k,w^k)=y^*-y^k\to 0,
$$
proving item (\ref{aw_gradLtnlp}). 

Let us prove item (\ref{aw_glbdg}). By the feasibility of $(x^*,y^*)$ we have 
$g_i(x^*)\leq 0$ for all $i=1,\ldots,m$. If $g_i(x^*)=0$, then 
$\min\{-g_i(x^k),\lambda_i^{g,k}\}\to 0$ since $g_i(x^k)\to 0$ and 
$\lambda_i^{g,k}\geq 0$. On the other hand, if $g_i(x^*)<0$, we may assume 
that $g_i(x^k)<0$ for all $k$. Thus, $g_i^+(x^k)=0$, yielding 
$\lambda_i^{g,k}=\rho_kg_i^+(x^k)=0$. Therefore, 
$\min\{-g_i(x^k),\lambda_i^{g,k}\}=0$. Items (\ref{aw_thlbdth}) and 
(\ref{aw_HtillbdHtil}) can be proved by the same reasoning. 

Now, note that by (\ref{aw_gradLpwG}) ,(\ref{aw_gradLpwH}) and (\ref{aw_gradphiw}) we have 
\begin{equation}
\label{aw_lbdGH}
\lambda^{G,k}=\mu^{0,k}*w^{H,k}\quad\mbox{and}\quad
\lambda^{H,k}=\mu^{H,k}-\mu^{0,k}*w^{G,k}.
\end{equation}
Therefore, using the fact that $w^k\in W$, we obtain 
$$
\lambda_i^{G,k}w_i^{G,k}=\mu_i^{0,k}w_i^{H,k}w_i^{G,k}=0
$$
for all $i=1,\ldots,m$. Furthermore, given $i\notin I_0(x^*)$, we have 
$$
w_i^{G,k}\to(w_i^*)^G=G_i(x^*)=x_i^*\neq 0,
$$ 
implying that $\lambda_i^{G,k}=0$ for all $k$ large enough. So,
$\min\{|G_i(x^k)|,|\lambda_i^{G,k}|\}=0$. On the other hand, if 
$i\in I_0(x^*)$, we have $G_i(x^k)\to G_i(x^*)=x_i^*=0$, and hence, 
$\min\{|G_i(x^k)|,|\lambda_i^{G,k}|\}\to 0$, proving item (\ref{aw_GlbdG}). 

To prove the next item, note that using (\ref{aw_lbdGH}), (\ref{aw_complw}) 
and the fact that $w^k\in W$,
\begin{equation}
\label{aw_lbdHwH}
\lambda^{H,k}*w^{H,k}=\mu^{H,k}*w^{H,k}-
\mu^{0,k}*w^{G,k}*w^{H,k}=0.
\end{equation}
By the feasibility of $(x^*,y^*)$, we have $H(y^*)\leq 0$. In the 
case $H_i(y^*)<0$, there holds 
$$
w_i^{H,k}\to(w_i^*)^H=-H_i(y^*)>0,
$$ 
giving $\lambda_i^{H,k}=0$ for all $k$ large enough. Thus, 
$\min\{-H_i(y^k),|\lambda_i^{H,k}|\}=0$. On the other hand, if 
$H_i(y^*)=0$, we have $H_i(y^k)\to H_i(y^*)=0$, and consequently, 
$\min\{-H_i(y^k),|\lambda_i^{H,k}|\}\to 0$, proving item (\ref{aw_HlbdH}) and  
completing the proof.
\endproof

\section{Relations to other sequential optimality conditions}
\label{aw_sec:rel_akkt_cakkt}
In this section we discuss the relationships between approximate 
stationarity for standard nonlinear optimization and $AW$-stationarity. 

As well known, every minimizer of an optimization problem is AKKT (see 
Definition \ref{aw_def:akkt}). However, and surprisingly, we start by proving 
that the AKKT condition fails to detect good candidates for optimality 
for every MPCaC problem.

\begin{theorem}
\label{aw_th:akktrelax}
Every feasible point $(\bar{x},\bar{y})$ for the relaxed problem (\ref{aw_prob:relax1})
is AKKT. 
\end{theorem}
\beginproof
We need to prove that there exist sequences $(x^k,y^k)\subset\R^n\times\R^n$ and 
$$
\big(\mu^{g,k},\mu^{h,k},\mu^{\theta,k},
\mu^{H,k},\mu^{\tilde{H},k},\mu^{\xi,k}\big)\subset\R_+^m\times\R^p\times\R_+
\times\R_+^n\times\R_+^n\times\R^n
$$
such that $(x^k,y^k)\to(\bar{x},\bar{y})$ and 
\begin{subequations}
\begin{align}
\left(\begin{array}{c} \nabla_xL(x^k,\mu^{g,k},\mu^{h,k}) \\ 0 \end{array}\right)+
\left(\begin{array}{c} 0 \\ \mu^{\theta,k}\nabla\theta(y^k) \end{array}\right)+
\sum_{i=1}^{n}\left(\begin{array}{c} 0 \\ \mu_i^{H,k} \nabla H_i(y^k)\end{array}\right) 
\notag{} \\ +\sum_{i=1}^{n}\left(\begin{array}{c} 0 \\ \mu_i^{\tilde{H},k}\nabla
\tilde{H}_i(y^k)\end{array}\right)+\sum_{i=1}^{n}\mu_i^{\xi,k}
\left(\begin{array}{c} H_i(y^k)\nabla G_i(x^k)\vspace{3pt} \\ 
G_i(x^k)\nabla H_i(y^k)\end{array}\right)\to 0, \label{aw_eq_aw_akktr1} \\ 
\min\{-g(x^k),\mu^{g,k}\}\to 0\,, \quad 
\min\{-\theta(y^k),\mu^{\theta,k}\}\to 0, \label{aw_eq_aw_akktr2} \\ 
\min\{-H(y^k),\mu^{H,k}\}\to 0\,,\quad
\min\{-\tilde{H}(y^k),\mu^{\tilde{H},k}\}\to 0.  \label{aw_eq_aw_akktr3} 
\end{align}
\end{subequations}
Let $b=\nabla f(\bar{x})$ and define $x^k=\bar{x}$, $\mu^{g,k}=0$, $\mu^{h,k}=0$, 
$\mu^{\theta,k}=0$, $\mu^{\tilde{H},k}=0$ and 
\begin{subequations}
\begin{align*}
y_i^k=\bar{y}_i\,,\,\, \mu_i^{H,k}=0\,,\,\, \mu_i^{\xi,k}=\dfrac{b_i}{y_i^k} 
\mbox{ for } i\in I_{0+}(\bar{x},\bar{y})\cup I_{01}(\bar{x},\bar{y}), \\ 
y_i^k=\dfrac{b_i}{k}\,,\,\, \mu_i^{H,k}=0\,,\,\, \mu_i^{\xi,k}=k 
\mbox{ for } i\in I_{00}(\bar{x},\bar{y}), \\ 
y_i^k=-\dfrac{{\rm sign}(\bar{x}_i)b_i}{k}\,,\,\,\mu_i^{\xi,k}=
-{\rm sign}(\bar{x}_i)k \,,\,\, \mu_i^{H,k}=-\mu_i^{\xi,k}x_i^k \mbox{ for } 
i\in I_{\pm 0}(\bar{x},\bar{y}).
\end{align*}
\end{subequations}
Thus we have $\mu_i^{H,k}\geq 0$, $(x^k,y^k)\to(\bar{x},\bar{y})$, 
$$
\nabla_{x_i}L(x^k,\mu^{g,k},\mu^{h,k})-\mu_i^{\xi,k}y_i^k=b_i-\mu_i^{\xi,k}y_i^k\to 0, 
$$
and 
$$
-\mu^{\theta,k}-\mu_i^{H,k}+\mu_i^{\tilde{H},k}-\mu_i^{\xi,k}x_i^k\to 0
$$
for all $i=1,\ldots,n$, giving (\ref{aw_eq_aw_akktr1}). Moreover, it is easy to see 
that (\ref{aw_eq_aw_akktr2}) and (\ref{aw_eq_aw_akktr3}) also hold. 
\endproof

Another sequential optimality condition for standard NLP is 
PAKKT (Definition \ref{aw_def:pakkt}). It is stronger than AKKT, but not 
stronger than $AW$-stationarity. The next example shows that PAKKT for the relaxed problem 
does not imply $AW$-stationarity, even under strict complementarity.

\begin{example}
\label{aw_ex:pakktnotaw}
Consider the MPCaC and the corresponding relaxed problem given below.
$$
\begin{array}{lr}
\begin{array}{cl}
\displaystyle\mathop{\rm minimize }_{x\in\R^2} & x_2  \\
{\rm subject\ to } & x_1^2 \leq 0, \\
& \|x\|_0\leq 1,\\
& \\
&
\end{array}
\hspace{.5cm}
&
\begin{array}{cl}
\displaystyle\mathop{\rm minimize }_{x,y\in\R^2} & x_2 \\
{\rm subject\ to } & x_1^2 \leq 0, \\
& y_1+y_2\geq 1, \\
& x_iy_i=0,\; i=1,2, \\
& 0\leq y_i \leq 1, \; i=1,2.
\end{array}
\end{array}
$$
Given $a>0$, we claim that the point $(\bar{x},\bar{y})$, with $\bar{x}=(0,a)$ and 
$\bar{y}=(1,0)$, is PAKKT but not $AW$-stationary. Indeed, for the first statement, 
consider the sequences $(x^k,y^k)\subset\R^2\times\R^2$ and 
$$
(\gamma^k)=\big(\lambda^{g,k},\lambda^{\theta,k},\mu^k,\lambda^{\tilde{H},k},
\lambda^{\xi,k}\big)\subset\R_+\times\R_+\times\R_+^2\times\R_+^2\times\R^2
$$
given by $x^k=(1/k^3,a)$, $y^k=(1,-1/k)$, $\lambda^{g,k}=k^2$, 
$\lambda^{\theta,k}=0$, $\mu^k=(0,ak)$, $\lambda^{\tilde{H},k}=(0,0)$ and 
$\lambda^{\xi,k}=(0,k)$. Then we have $(x^k,y^k)\to(\bar{x},\bar{y})$ and, 
denoting $\xi(x,y)=x*y$, the gradient of the Lagrangian of the relaxed problem 
reduces to 
\begin{align*}
\left(\begin{array}{c}\nabla f(x^k) \\ 0 \end{array}\right)
+\lambda^{g,k}
\left(\begin{array}{c}\nabla g(x^k) \\ 0 \end{array}\right)
+\mu_2^k \nabla H_2(y^k)+\lambda_2^{\xi,k}\nabla\xi_2(x^k,y^k) \\
=\left(\begin{array}{c} 0 \\ 1 \\ 0 \\ 0 \end{array}\right)+
\left(\begin{array}{c} 2\lambda^{g,k}x_1^k \\ 0 \\ 0 \\ 0 
\end{array}\right)+
\left(\begin{array}{c} 0 \\ 0 \\ 0 \\ -\mu_2^k \end{array}
\right)+\left(\begin{array}{c} 0 \\ 
\lambda_2^{\xi,k}y_2^k \\ 0 \\ \lambda_2^{\xi,k}x_2^k \end{array}\right)
=\left(\begin{array}{c} 2/k \\ 0 \\ 0 \\ 0 
\end{array}\right)\to 0, 
\end{align*}
proving (\ref{aw_gradpakkt}). Now, note that 
$g(x^k)\to g(\bar{x})=0$ and $\theta(y^k)\to\theta(\bar{y})=0$, 
which in turn imply that 
\begin{equation}
\label{aw_eq_ex:akktnotaw1}
\min\{-g(x^k),\lambda^{g,k}\}\to 0\quad\mbox{and}\quad
\min\{-\theta(y^k),\lambda^{\theta,k}\}\to 0.
\end{equation}
Moreover, we have $-\tilde{H}(y^k)\to -\tilde{H}(\bar{y})\geq 0$ and 
$\lambda^{\tilde{H},k}=(0,0)$, giving 
\begin{equation}
\label{aw_eq_ex:akktnotaw2}
\min\{-\tilde{H}(y^k),\lambda^{\tilde{H},k}\}\to 0.
\end{equation}
Furthermore, since $-H_1(y^k)\to -H_1(\bar{y})\geq 0$, $\mu_1^k=0$ and 
$-H_2(y^k)=y_2^k\to 0$, we have 
\begin{equation}
\label{aw_eq_ex:akktnotaw3}
\min\{-H(y^k),\mu^k\}\to 0.
\end{equation}
Conditions (\ref{aw_eq_ex:akktnotaw1}), (\ref{aw_eq_ex:akktnotaw2}) and 
(\ref{aw_eq_ex:akktnotaw3}) prove the approximate complementarity (\ref{aw_complpakkt}). 
Moreover, we have $\delta_k= \|(1,\gamma^k)\|_\infty=k^2$ for 
all $k$ large enough, 
$$
\displaystyle\mathop{\rm lim\, sup}_{k\to\infty}\frac{\lambda^{g,k}}{\delta_k}>0
\quad\mbox{and}\quad \lambda^{g,k}g(x^k)>0.
$$
For the remaining multipliers the ${\rm lim\, sup}$ is zero and so we conclude that 
(\ref{aw_pospakkt1}) and (\ref{aw_pospakkt2}) hold, proving that Definition \ref{aw_def:pakkt} 
is satisfied, that is, $(\bar{x},\bar{y})$ is PAKKT. 

Now, let us see that $(\bar{x},\bar{y})$ is not $AW$-stationary. For this purpose, assume 
that the sequences $(x^k,y^k)\subset\R^2\times\R^2$ and 
$$
(\lambda^k)=\big(\lambda^{g,k},\lambda^{\theta,k},
\lambda^{G,k},\lambda^{H,k},\lambda^{\tilde{H},k}\big)\subset\R_+\times\R_+
\times\R^2\times\R^2\times\R_+^2
$$
are such that $(x^k,y^k)\to(\bar{x},\bar{y})$ and 
$\min\{|G_2(x^k)|,|\lambda_2^{G,k}|\}\to 0$. Then, since $$|G_2(x^k)|=|x_2^k|\to a>0,$$ 
we obtain $\lambda_2^{G,k}\to 0$. Therefore, the expression 
$$
\nabla_xL(x^k,\lambda^{g,k})+\sum_{i=1}^2\lambda_i^{G,k}
\nabla G_i(x^k)=\left(\begin{array}{c} 2\lambda^{g,k}x_1^k+
\lambda_1^{G,k} \\ 1+\lambda_2^{G,k} \end{array}\right)
$$
cannot converge to zero. Thus, taking into account (\ref{aw_nablaL0}), 
item (\ref{aw_gradLtnlp}) of Definition~\ref{aw_def:aw} does not hold and 
hence $(\bar{x},\bar{y})$ is not $AW$-stationary. 
\end{example}

\medskip

In contrast to AKKT and PAKKT, the other classical sequential optimality 
condition, CAKKT (Definition \ref{aw_def:cakkt}), does imply $AW$-stationarity, 
as we can see in the next result. 

\begin{theorem}
\label{aw_th:cakkt_aw}
If $(\bar{x},\bar{y})$ is a CAKKT point for the relaxed problem (\ref{aw_prob:relax1}), 
then it is $AW$-stationary. 
\end{theorem}
\beginproof
In view of Definition \ref{aw_def:cakkt}, there exist 
sequences $(x^k,y^k)\subset\R^n\times\R^n$ and 
$$
\big(\lambda^{g,k},\lambda^{h,k},\lambda^{\theta,k},
\mu^k,\lambda^{\tilde{H},k},\lambda^{\xi,k}\big)\subset\R_+^m\times\R^p\times\R_+
\times\R_+^n\times\R_+^n\times\R^n
$$
such that $(x^k,y^k)\to(\bar{x},\bar{y})$, 
\begin{subequations}
\begin{align}
\left(\begin{array}{c}\nabla_xL(x^k,\lambda^{g,k},\lambda^{h,k}) \\ 0 \end{array}\right)+
\left(\begin{array}{c} 0 \\ \lambda^{\theta,k}\nabla\theta(y^k) \end{array}\right)+
\sum_{i=1}^{n}\left(\begin{array}{c} 0 \\ \mu_i^k \nabla H_i(y^k)\end{array}\right) 
\notag{} \\ 
+\sum_{i=1}^{n}\left(\begin{array}{c} 0 \\ \lambda_i^{\tilde{H},k}\nabla
\tilde{H}_i(y^k)\end{array}\right)+\sum_{i=1}^{n}\lambda_i^{\xi,k}
\left(\begin{array}{c} H_i(y^k)\nabla G_i(x^k)\vspace{3pt} \\ 
G_i(x^k)\nabla H_i(y^k)\end{array}\right)\to 0, \label{aw_eq_cakkt_aw1} \\ 
\lambda^{g,k}*g(x^k)\to 0\,, \quad \lambda^{h,k}*h(x^k)\to 0\,, \quad 
\lambda^{\theta,k}\theta(y^k)\to 0\,, \label{aw_eq_cakkt_aw2} \\ 
\mu^k*H(y^k)\to 0\,,\quad \lambda^{\tilde{H},k}*\tilde{H}(y^k)\to 0, \label{aw_eq_cakkt_aw3} \\
\lambda^{\xi,k}*G(x^k)*H(y^k)\to 0. \label{aw_eq_cakkt_aw4} 
\end{align}
\end{subequations}
So, we may define 
$$
\lambda^{H,k}=\mu^k+\lambda^{\xi,k}*G(x^k)\quad\mbox{and}
\quad\lambda^{G,k}=\lambda^{\xi,k}*H(y^k)
$$
to obtain item (\ref{aw_gradLtnlp}) of Definition \ref{aw_def:aw} from (\ref{aw_eq_cakkt_aw1}). 
Items (\ref{aw_glbdg}), (\ref{aw_thlbdth}) and (\ref{aw_HtillbdHtil}) follow 
from (\ref{aw_eq_cakkt_aw2}), (\ref{aw_eq_cakkt_aw3}) and Remark \ref{aw_rm:compl4_compl1}. 
Let us prove item (\ref{aw_GlbdG}). For $i\in I_{0}$, there holds 
$$G_i(x^k)\to G_i(\bar{x})=\bar{x}_i=0.$$ Thus, 
$\min\{|G_i(x^k)|,|\lambda_i^{G,k}|\}\to 0$. If $i\notin I_{0}$, we have 
$G_i(x^k)\to\bar{x}_i\neq 0$, which in view of (\ref{aw_eq_cakkt_aw4}) yields 
$$
\lambda_i^{G,k}=\lambda_i^{\xi,k}H_i(y^k)\to 0.
$$
Therefore, $\min\{|G_i(x^k)|,|\lambda_i^{G,k}|\}\to 0$ for all $i=1,\ldots,n$. 
Finally, in order to prove item (\ref{aw_HlbdH}), note that 
(\ref{aw_eq_cakkt_aw3}) and (\ref{aw_eq_cakkt_aw4}) give  
$$
\lambda_i^{H,k}H_i(y^k)=\mu_i^kH_i(y^k)+\lambda_i^{\xi,k}G_i(x^k)H_i(y^k)\to 0.
$$
So, applying the argument of Remark \ref{aw_rm:compl4_compl1} with 
$\alpha^k=|\lambda_i^{H,k}|$ and $\beta^k=H_i(y^k)$, we obtain 
$$
\min\{-H_i(y^k),|\lambda_i^{H,k}|\}\to 0
$$
for all $i=1,\ldots,n$. Therefore, $(\bar{x},\bar{y})$ is 
$AW$-stationary for the problem (\ref{aw_prob:relax1}).
\endproof

\begin{remark}
\label{aw_rm:cakkt}
Despite being stronger, we emphasize that the sequential 
optimality condition CAKKT is not so suitable to deal with MPCaC problems as 
$AW$-stationarity. The goal of considering CAKKT is to obtain, under certain 
constraint qualifications, KKT points for standard nonlinear programming problems. 
However, as we have been discussed, MPCaC are very degenerate problems because of 
the problematic complementarity constraint $G(x)*H(y)=0$. This means that we cannot 
expect to find strong stationary points for this class of problems and thereby 
making $AW$-stationarity a good tool for dealing with them. 
\end{remark}

To finish this section, we relate our sequential optimality condition to the 
tightened problem. The following result is a sequential version of 
Proposition~\ref{aw_prop:wstat_kkttnlp}. 

\begin{theorem}
\label{aw_th:aw_akkt}
Let $(\bar{x},\bar{y})$ be a feasible point of the relaxed problem (\ref{aw_prob:relax1}). 
Then $(\bar{x},\bar{y})$ is $AW$-stationary if and only if it is an AKKT point for the 
tightened problem TNLP$_{I_0}(\bar{x},\bar{y})$ defined in (\ref{aw_prob:tight}).
\end{theorem}
\beginproof
Suppose first that $(\bar{x},\bar{y})$ is $AW$-stationary. Then, in view of 
Lemma \ref{aw_lm:awI0}, we conclude that there exist sequences 
$(x^k,y^k)\subset\R^n\times\R^n$ and 
$$
(\lambda^k)=\big(\lambda^{g,k},\lambda^{h,k},\lambda^{\theta,k},
\lambda^{G,k},\lambda^{H,k},\lambda^{\tilde{H},k}\big)
\subset\R_+^m\times\R^p\times\R_+\times\R^n\times\R^n\times\R_+^n
$$
such that $(x^k,y^k)\to(\bar{x},\bar{y})$, 
\begin{subequations}
\begin{align}
\nabla_xL(x^k,\lambda^{g,k},\lambda^{h,k})+\sum_{i\in I_0}
\lambda_i^{G,k}\nabla G_i(x^k)\to 0, \label{aw_eq_aw_akkt1} \\ 
\lambda^{\theta,k}\nabla\theta(y^k) +\sum_{i=1}^{n}\lambda_i^{H,k} 
\nabla H_i(y^k)+\sum_{i=1}^{n}\lambda_i^{\tilde{H},k}\nabla\tilde{H}_i(y^k)\to 0, 
\label{aw_eq_aw_akkt2} \\ 
\min\{-g(x^k),\lambda^{g,k}\}\to 0\,, \quad \min\{-\theta(y^k),\lambda^{\theta,k}\}\to 0,
\label{aw_eq_aw_akkt3} \\ \min\{-H_i(y^k),|\lambda_i^{H,k}|\}\to 0\,,\ i=1,\ldots,n,
\quad\min\{-\tilde{H}(y^k),\lambda^{\tilde{H},k}\}\to 0.  \label{aw_eq_aw_akkt4} 
\end{align}
\end{subequations}
For $i\in I_{0+}\cup I_{01}$ we have 
$
H_i(y^k)\to H_i(\bar{y})=-\bar{y}_i<0.
$
Therefore, we can assume without loss of generality that there exists $\epsilon>0$ 
such that $-H_i(y^k)\geq\epsilon$ for all $k$. So, using (\ref{aw_eq_aw_akkt4}), we 
obtain $|\lambda_i^{H,k}|\to 0$, which in turn implies that 
$$
\sum_{i\in I_{0+}\cup I_{01}}\lambda_i^{H,k}\nabla H_i(y^k)\to 0.
$$
By subtracting this from (\ref{aw_eq_aw_akkt2}), we obtain 
$$
\lambda^{\theta,k}\nabla\theta(y^k) +\sum_{i\in I_{00}\cup I_{\pm 0}}
\lambda_i^{H,k} \nabla H_i(y^k)+\sum_{i=1}^{n}
\lambda_i^{\tilde{H},k}\nabla\tilde{H}_i(y^k)\to 0.
$$
So, we can redefine $\lambda_i^{H,k}$, $i\in I_{0+}\cup I_{01}$, 
to be zero, without affecting (\ref{aw_eq_aw_akkt2}). 
Therefore, taking into account (\ref{aw_eq_aw_akkt1}), (\ref{aw_eq_aw_akkt3}), the second 
part of (\ref{aw_eq_aw_akkt4}) and the fact that 
$\min\{-H_i(y^k),\lambda_i^{H,k}\}=0$ for $i\in I_{0+}\cup I_{01}$, 
we conclude that $(\bar{x},\bar{y})$ is AKKT for TNLP$_{I_0}(\bar{x},\bar{y})$, 
which we recall here for convenience, 
$$
\begin{array}{cl}
\displaystyle\mathop{\rm minimize }_{x,y} & f(x) \\
{\rm subject\ to } & g(x)\leq 0, h(x)=0, \\
& \theta(y)\leq 0, \\
& G_i(x)=0,\; i \in I_{0}, \\
& H_i(y)\leq 0,\; i\in I_{0+}\cup I_{01}, \\
& H_i(y)=0,\; i\in I_{00}\cup I_{\pm0}, \\
& \tilde{H}(y)\leq 0.
\end{array}
$$

To prove the converse, suppose that $(\bar{x},\bar{y})$ is AKKT for 
TNLP$_{I_0}(\bar{x},\bar{y})$. Then there exist sequences $(x^k,y^k)\subset\R^n\times\R^n$ 
and 
$$
(\lambda^k)=\big(\lambda^{g,k},\lambda^{h,k},\lambda^{\theta,k},
\lambda_{I_{0}}^{G,k},\lambda^{H,k},\lambda^{\tilde{H},k}\big)
\subset\R_+^m\times\R^p\times\R_+\times\R^{|I_{0}|}\times\R^n\times\R_+^n,
$$
with $\lambda_i^{H,k}\geq 0$ for $i\in I_{0+}\cup I_{01}$, such that 
$(x^k,y^k)\to(\bar{x},\bar{y})$, 
\begin{subequations}
\begin{align}
\nabla_xL(x^k,\lambda^{g,k},\lambda^{h,k})+\sum_{i\in I_0}
\lambda_i^{G,k}\nabla G_i(x^k)\to 0, \label{aw_eq_akkt_aw1} \\ 
\lambda^{\theta,k}\nabla\theta(y^k) +\sum_{i=1}^{n}\lambda_i^{H,k} 
\nabla H_i(y^k)+\sum_{i=1}^{n}\lambda_i^{\tilde{H},k}\nabla\tilde{H}_i(y^k)\to 0, 
\label{aw_eq_akkt_aw2} \\ 
\min\{-g(x^k),\lambda^{g,k}\}\to 0\,, \quad \min\{-\theta(y^k),
\lambda^{\theta,k}\}\to 0, \label{aw_eq_akkt_aw3} \\ 
\min\{-H_i(y^k),\lambda_i^{H,k}\}\to 0,\ i\in I_{0+}\cup I_{01}\,,
\quad\min\{-\tilde{H}(y^k),\lambda^{\tilde{H},k}\}\to 0. \label{aw_eq_akkt_aw4} 
\end{align}
\end{subequations}
Extending the sequence $\big(\lambda_{I_{0}}^{G,k}\big)$ from $\R^{|I_{0}|}$ to 
$\R^n$ by letting $\lambda_i^{G,k}=0$ for $i\notin I_{0}$, we can 
rewrite (\ref{aw_eq_akkt_aw1}) as 
\begin{align}
\nabla_xL(x^k,\lambda^{g,k},\lambda^{h,k})+\sum_{i=1}^{n}
\lambda_i^{G,k}\nabla G_i(x^k)\to 0. \label{aw_eq_akkt_aw5}
\end{align}
Moreover, for $i\in I_{0}$, there holds $G_i(x^k)\to G_i(\bar{x})=\bar{x}_i=0$. 
Thus, 
\begin{equation}
\label{aw_eq_akkt_aw6}
\min\{|G_i(x^k)|,|\lambda_i^{G,k}|\}\to 0
\end{equation}
for all $i=1,\ldots,n$. Besides, for $i\in I_{00}\cup I_{\pm 0}$, we have 
$
H_i(y^k)\to H_i(\bar{y})=-\bar{y}_i=0,
$
which implies $\min\{-H_i(y^k),|\lambda_i^{H,k}|\}\to 0$. Therefore, in view of 
(\ref{aw_eq_akkt_aw4}) and the fact that $\lambda_i^{H,k}\geq 0$ for 
$i\in I_{0+}\cup I_{01}$, we have 
\begin{equation}
\label{aw_eq_akkt_aw7}
\min\{-H_i(y^k),|\lambda_i^{H,k}|\}\to 0
\end{equation}
for all $i=1,\ldots,n$. Thus, from (\ref{aw_eq_akkt_aw2}), (\ref{aw_eq_akkt_aw3}), the 
second part of (\ref{aw_eq_akkt_aw4}), (\ref{aw_eq_akkt_aw5}), (\ref{aw_eq_akkt_aw6}) and 
(\ref{aw_eq_akkt_aw7}), we conclude that 
$(\bar{x},\bar{y})$ satisfies the conditions of Definition \ref{aw_def:aw}, that is, 
$(\bar{x},\bar{y})$ is an $AW$-stationary point for the problem (\ref{aw_prob:relax1}).
\endproof

\section{Conclusion}
\label{aw_sec:concl} 
In this paper we have presented a sequential optimality condition, namely Approximate Weak stationarity ($AW$-stationarity), for Mathematical Programs with Cardinality Constraints (MPCaC). 
This condition improves $W_I$-statio\-na\-rity, which was established in our previous work~\cite{KrulikovskiRibeiroSachine20aX}. 

Several theoretical results were presented, such as: $AW$-stationarity is a legitimate 
optimality condition independently of any constraint qualification; every feasible 
point of MPCaC is AKKT; the equivalence between the $AW$-stationarity and AKKT for 
the tightened problem TNLP$_{I_0}$.
In addition, we have established some relationships between our $AW$-stationarity 
and other usual sequential optimality conditions, such as AKKT, CAKKT and PAKKT, 
by means of properties, examples and counterexamples. 

It should be mentioned that, even though the computational appeal of the sequential 
optimality conditions, in this work we were not concerned with algorithmic consequences, 
which is subject of ongoing research. 

%
%
%
%



\end{document}